\documentclass{amsart}%
\usepackage{amsfonts}
\usepackage{amsmath}
\usepackage{amssymb}
\usepackage{graphicx}%
\newtheorem{theorem}{Theorem}
\theoremstyle{plain}

\newtheorem{conjecture}{Conjecture}

\newtheorem{example}{Example}
\newtheorem{lemma}{Lemma}

\newtheorem{proposition}{Proposition}
\newtheorem{remark}{Remark}

\theoremstyle{definition}
\newtheorem*{acknowledgement}{Acknowledgments}

\theoremstyle{remark}

\numberwithin{equation}{section}
\begin{document}
\title[Expanding homogeneous Ricci solitons]{Linear stability of homogeneous Ricci solitons}
\author{Christine Guenther}
\address[Christine Guenther]{ Pacific University}
\email{guenthec@pacificu.edu}
\author{James Isenberg}
\address[James Isenberg]{ University of Oregon}
\email{jim@newton.uoregon.edu}
\urladdr{http://physics.uoregon.edu/\symbol{126}jim/}
\author{Dan Knopf}
\address[Dan Knopf]{ University of Texas at Austin}
\email{danknopf@math.utexas.edu}
\urladdr{http://www.ma.utexas.edu/\symbol{126}danknopf/}
\thanks{Revised September 5, 2006}

\begin{abstract}
As a step toward understanding the analytic behavior of Type-III Ricci flow
singularities, i.e.~immortal solutions that exhibit $|\operatorname*{Rm}|\leq
C/t$ curvature decay, we examine the linearization of an equivalent flow at
fixed points discovered recently by Baird--Danielo and Lott: nongradient
homogeneous expanding Ricci solitons on nilpotent or solvable Lie groups. For
all explicitly known nonproduct examples, we demonstrate linear stability of
the flow at these fixed points and prove that the linearizations generate
$C_{0}$ semigroups.

\end{abstract}
\maketitle

\section{Introduction}

In \cite{GIK02}, the authors investigate the stability of compact flat and
Ricci-flat solutions of Ricci flow --- where there is a center manifold
present in the space of Riemannian metrics --- by a method incorporating the
linearized Ricci flow operator, the analytic semigroup it generates, and
appropriate interpolation spaces. (One of our results has recently been
strengthened by \v{S}e\v{s}um using somewhat different techniques; see
\cite{Sesum04}.)

The work described in this note concerns a related question: stability along a
collapsing locally homogeneous solution of Ricci flow. One incentive for
considering this more nuanced issue arises when studying the geometrization of
a closed $3$-manifold. The idea is essentially as follows. (See
Section~\ref{Motivation} for a more precise discussion.) Let $\mathcal{M}%
_{t}^{3}$ denote the (possibly disconnected) solution of Ricci
flow-with-surgery at some non-surgery time $t\gg0$. A central claim of the
Hamilton--Perelman Geometrization program is that there exists a decomposition
$\mathcal{M}_{t}^{3}=\mathcal{M}_{\mathrm{thick}}^{3}\cup\mathcal{M}%
_{\mathrm{thin}}^{3}$. (See \cite{Perelman1, Perelman2}, \cite{KL}, and
\cite{CZ06}.) $\mathcal{M}_{\mathrm{thick}}^{3}$ consists of those components
where the injectivity radius is controlled in a suitable manner; these
components should be close in a precise sense to a finite collection of
complete finite-volume hyperbolic $3$-manifolds, truncated along cuspidal tori
whose images are incompressible in $\mathcal{M}_{t}^{3}$. (Compare
\cite{Hamilton99}.) The components of $\mathcal{M}_{\mathrm{thin}}^{3}$ are
collapsed with local lower bounds on sectional curvature; these components
should be Cheeger--Gromov graph manifolds. (See \cite{SY00, SY05}.) A graph
manifold in the sense of Cheeger--Gromov \cite{CG86, CG90} is a closed
$3$-manifold admitting a decomposition by not necessarily incompressible tori
such that each complementary piece is a Seifert space. Any such manifold
either is a topological graph manifold (i.e.~it may be decomposed along
incompressible tori such that each complementary piece is a Seifert manifold)
or else is a connected sum of topological graph manifolds with lens space and
$\mathcal{S}^{1}\times\mathcal{S}^{2}$ factors.

Many intriguing questions remain open about the behavior of Ricci flow on
$\mathcal{M}_{\mathrm{thin}}^{3}$. A long-standing question of Hamilton asks
whether the Ricci flow on these components may be modeled in a suitable manner
by its behavior on the homogeneous Thurston model geometries. If so, one would
expect that the Ricci flow of the locally homogeneous Thurston geometries ---
among which $\operatorname*{nil}$ and $\operatorname*{sol}$ exhibit
Cheeger--Gromov collapse --- should in an appropriate sense be stable in the
space of $3$-dimensional solutions. (See \cite{IJ92} and \cite{KM01}, as well
as \cite{HI93} and \cite{Knopf00}.)

Recent work of Lott \cite{Lott06} allows a refinement of Hamilton's question.
In any dimension, a Type-III solution of Ricci flow $(\mathcal{M}^{n},g(t))$
is one that exists for $t\in\lbrack0,\infty)$ and satisfies $\sup
_{\mathcal{M}^{n}\times\lbrack0,\infty)}t|\operatorname*{Rm}|<\infty$. Given
such a solution, fix an origin $x\in\mathcal{M}^{n}$. Then for each $s>0$,
there is a rescaled pointed solution of Ricci flow $(\mathcal{M}^{n}%
,g_{s}(t),x)$ defined for all $t\in\lbrack0,\infty)$ by $g_{s}(t):=s^{-1}%
g(st)$. If the original solution is a finite-volume locally homogeneous
$3$-manifold, Lott proves that $\lim_{s\rightarrow\infty}g_{s}(\cdot)$ exists
as a homogeneous Ricci soliton on a $3$-dimensional \'{e}tale
groupoid.\footnote{More generally, Lott's compactness theorem
\cite[Theorem~1.4]{Lott06} shows that any sequence $(\mathcal{M}_{k}^{n}%
,g_{k}(\cdot),x_{k})$ of complete pointed Ricci flows defined on a common time
interval and obeying uniform curvature bounds subconverges to a solution of
Ricci flow on a pointed $n$-dimensional \'{e}tale groupoid. In general, one
does not know this is a soliton.} (Compare \cite{Glickenstein03}.) Let
$\tilde{g}(\cdot)$ denote the lift of $g(\cdot)$ to the universal cover
$\mathcal{\tilde{M}}^{3}$. The limit solution on $\mathcal{\tilde{M}}_{\infty
}^{3}$ is $\tilde{g}_{\infty}(t):=\lim_{s\rightarrow\infty}(\psi_{s}^{\ast
}(s^{-1}\tilde{g}(st)))$, where $\{\psi_{s}\}$ is a family of time-independent
diffeomorphisms. Thus $\tilde{g}_{\infty}(\cdot)$ encodes the asymptotic
geometry of $\tilde{g}(\cdot)$ at large times and large length
scales.\footnote{The limit $\tilde{g}_{\infty}(\cdot)$ should be compared to
the asymptotic gradient shrinking soliton of an ancient $\kappa$-solution,
which Perelman constructs \cite[\S 11.2]{Perelman1} as a blow-down limit as
$t\rightarrow-\infty$.} Lauret \cite{Lauret01} has shown that, in many cases,
a homogeneous soliton is essentially unique.\footnote{See
Theorem~\ref{Lauret2}, below.} Any such soliton is a fixed point of the
transformation $\tilde{g}(\cdot)\mapsto\tilde{g}_{\infty}(\cdot)$. Therefore,
one implication of Lott's results may be interpreted as the statement that if
a fixed $3$-dimensional Lie group admits a homogeneous Ricci soliton --- which
is true for $\operatorname*{nil}$ and $\operatorname*{sol}$, but not for
$\widetilde{\operatorname*{SL}}(2,\mathbb{R)}$ --- then the soliton acts in
this sense as an infinite-time attractor among all solutions respecting the
group symmetry. From the point of view of Hamilton's original question, one
may further ask whether a given homogeneous Ricci soliton acts as an
infinite-time attractor among \emph{all} nearby solutions, symmetric or not.
If so, then the collapsing behavior of the soliton is stable for Ricci flow,
in a precise sense.

For example, suppose $(\mathcal{M}^{3},g(t))$ is a locally homogeneous
$\operatorname*{nil}$-geometry solution of Ricci flow on a compact manifold.
Lifting to the universal cover $\mathcal{\tilde{M}}^{3}\approx\mathbb{R}^{3}$,
one finds homogeneous $\operatorname*{nil}$-geometry metrics $(\mathcal{\tilde
{M}}^{3},G(t))$ evolving by Ricci flow. In the six-dimensional
family\thinspace\footnote{A homogeneous solution of Ricci flow is determined
by its initial datum, a left-invariant metric. A left-invariant metric on an
$n$-dimensional homogeneous space is determined by its value on any tangent
space, hence may be identified with a positive-definite symmetric $(n\times
n)$-matrix. Any such matrix is determined by $n(n+1)/2$ parameters.}
$\mathcal{G}$ parameterizing such evolutions, exactly one is a
soliton.\footnote{Uniqueness was proved by Lauret \cite{Lauret01}. Explicit
constructions were done independently by Baird--Danielo \cite{BD06} and Lott
\cite{Lott06}. See Sections~\ref{Motivation} and \ref{Linear} for a more
thorough exposition.} Call this $G_{0}(\cdot)$, the lift of $g_{0}(\cdot)$.
(Note that $G_{0}(\cdot)$ is unique up to scaling, and that $g_{0}(\cdot)$ is
a locally homogeneous $\operatorname*{nil}$-geometry solution.) Lott shows
that the solution determined by the initial datum $G_{0}(0)$ is stable within
$\mathcal{G}$, in a precise sense. Is there anything one can say about the
stability of $G_{0}(0)$ among its non-homogeneous perturbations $\tilde{G}$,
or of $g_{0}(0)$ among its non-locally-homogeneous perturbations $\tilde{g}$?

In this note, we explore an approach to these questions from the point of view
of \cite{GIK02}. A homogeneous Ricci soliton may be regarded as a fixed point
of a related flow, solutions of which are in one-to-one correspondence with
solutions of Ricci flow. If one could show that the linearization of this flow
obeys a spectral bound and generates an analytic semigroup, then results of Da
Prato and Lunardi \cite[Theorem~2.2]{DPL88} would imply asymptotic stability
of the soliton among all nearby solutions of Ricci flow in a suitable
interpolation space. However, for the solitons under consideration, the
analytic and technical challenges of this program are formidable. The
linearizations of the modified flow at these solitons are second-order
tensor-valued differential operators with quadratically unbounded first and
zeroth order coefficients posed on noncompact manifolds. These operators are
not self adjoint in any weighted $L^{2}$ space. Such operators have
substantial independent interest from the perspective of analysis. Indeed, the
analysis of scalar-valued operators with unbounded coefficients is a very
active area of current research, because of their connection with
Ornstein--Uhlenbeck operators with unbounded drift and Schr\"{o}dinger
operators with unbounded potentials. For a small sampling, see \cite{CF04},
\cite{CV87}, \cite{GMV02}, \cite{LMP05}, \cite{MPSR05}, and \cite{Rabier05}.

From the point of view of analysis, the key results of this paper concern
certain geometrically motivated second-order non-self-adjoint tensor-valued
differential operators with quadratically unbounded lower-order coefficients
on noncompact manifolds. We show that these operators generate strongly
continuous semigroups and that their spectra have strictly negative
(i.e.~stable) real parts.

This note is organized as follows. In Section~\ref{Motivation}, we provide a
more thorough background and motivation for studying the asymptotic stability
of expanding homogeneous solitons. In Section~\ref{Linear}, we derive the
modified Ricci flow (\ref{MRF}) and its linearization (\ref{L}) at soliton
fixed points, with an eye toward studying this stability. In
Section~\ref{Spectra}, we discuss all three explicitly known nonproduct
examples: two are nilsolitons and one exhibits $\operatorname*{sol}$ geometry.
We prove linear stability of the modified Ricci flow at these fixed points.
(See Propositions~\ref{Stable3}, \ref{StableSol}, and \ref{Stable4}.) This
provides a formal argument in support of asymptotic stability. In
Section~\ref{Semigroups}, we prove that the linear operators generate $C_{0}$
semigroups. (See Theorem~\ref{StronglyContinuous}.) Even this result is far
from trivial and may be of independent analytic interest, as we have explained
in the previous two paragraphs. Our progress in this note falls short of
establishing analytic semigroups. So the question of asymptotic stability
remains open but appears promising, in light of our results in
Sections~\ref{Spectra} and \ref{Semigroups}. We hope that our conclusions here
facilitate further study of this interesting area of geometric analysis. In
Section~\ref{Open}, we briefly discuss some directions for further research.

\begin{acknowledgement}
C.G.~is partially supported by the Thomas J.~and Joyce Holce Professorship in
Science. J.I.~is partially supported by NSF grant PHY-0354659. D.K.~is
partially supported by NSF grants DMS-0511184 and DMS-0505920.

D.K.~ thanks Sigurd Angenent and Rafe Mazzeo for their interest and discussions.

The authors thank the referee for suggestions to improve the exposition.
\end{acknowledgement}

\section{Motivation and background\label{Motivation}}

\subsection{Ricci soliton structures}

A Ricci soliton structure consists of the data $(\mathcal{M}^{n},g,X,\alpha)$,
where $\mathcal{M}^{n}$ is a smooth connected manifold, $g$ is a complete
Riemannian metric on $\mathcal{M}^{n}$, $X$ is a complete vector field on
$\mathcal{M}^{n}$, and $\alpha$ is a constant such that%
\begin{equation}
2\operatorname*{Rc}+\alpha g+\mathcal{L}_{X}g=0. \label{RSSE}%
\end{equation}
One says the soliton is shrinking, steady, or expanding if $\alpha<0$,
$\alpha=0$, or $\alpha>0$, respectively. If $\eta_{t}$ is the one-parameter
family of diffeomorphisms generated by the vector fields $Y(\cdot,t)=(1+\alpha
t)^{-1}X(\cdot)$, then $\bar{g}(\cdot,t):=(1+\alpha t)\eta_{t}^{\ast}g(\cdot)$
is a solution of Ricci flow that evolves self similarly, i.e.~only by scaling
and diffeomorphism. Conversely, after normalization, every self similar
solution of Ricci flow arises in this manner.

\subsection{Solitons as finite-time singularity models}

A major reason for studying shrinking or steady solitons is the information
they can provide about finite-time singularities. If $(\mathcal{M}^{n},g(t))$
is a compact solution of Ricci flow that exists up to a maximal time
$T<\infty$, one may study the developing singularity by forming a sequence of
pointed dilated solutions $(\mathcal{M}^{n},g_{k}(t),x_{k})$, where
$g_{k}(t)=\lambda_{k}g(t_{k}+t/\lambda_{k})$ for some sequence of positive
constants $\lambda_{k}\rightarrow\infty$. To get a complete, smooth limit, one
takes $t_{k}\nearrow T$ and chooses the dilation factors $\lambda_{k}$ to be
comparable to the supremum of $\left\vert \operatorname*{Rm}\right\vert $ over
appropriate space-time neighborhoods of $(x_{k},t_{k})$. Under quite general
conditions, the injectivity radius estimate given by Perelman's No Local
Collapsing Theorem \cite{Perelman1} then allows application of Hamilton's
Compactness Theorem \cite{Hamilton95a} to show that the sequence
$(\mathcal{M}^{n},g_{k}(t),x_{k})$ subconverges, modulo time-independent
diffeomorphisms $\psi_{k}$, to a pointed limit solution of Ricci flow
$(\mathcal{M}_{\infty}^{n},g_{\infty}(t),x_{\infty})$, called a singularity
model. In dimension three, every finite-time singularity admits gradient
solitons among its possible singularity models. (See \cite{Hamilton95} and
\cite{Perelman1}.) In any dimension, it is reasonable to expect finite-time
singularity models to be shrinking or steady solitons. For example, the
neckpinch singularity \cite{AK1, AK2} is modeled by the shrinking gradient
soliton $(\mathbb{R}\times S^{n},ds^{2}+2(n-1)g_{\operatorname*{can}%
},\operatorname*{grad}(s^{2}/4),-1)$ for all $n\geq2$.

\subsection{Solitons as infinite-time singularity models}

Our interest in this note is in immortal solutions of Ricci flow. As noted
above, a Type-III solution of Ricci flow $(\mathcal{M}^{n},g(t))$ is an
immortal solution that satisfies $\sup_{\mathcal{M}^{n}\times\lbrack0,\infty
)}t|\operatorname*{Rm}|<\infty$. There are many interesting open questions
about the analytic properties of Type-III solutions in all dimensions $n\geq3$.

For Type-III solutions, no uniform injectivity radius estimate is known. In
fact, none is possible. For example, compact locally homogeneous manifolds
with $\operatorname*{nil}^{3}$ geometry occur as mapping tori of $\Upsilon
_{k}:T^{2}\rightarrow T^{2}$ induced by $%
\begin{pmatrix}
1 & k\\
0 & 1
\end{pmatrix}
\in\operatorname*{SL}(2,\mathbb{Z})$ with $k\neq0$. As shown in \cite{IJ92}
and \cite{KM01}, Ricci flow solutions with this geometry are all Type-III,
with $|\operatorname*{Rm}(g(t))|\leq C/(1+t)$ and with injectivity radius
decaying like $t^{-1/6}$ as $t\rightarrow\infty$. Other three-dimensional
examples are analyzed in \cite{IJ92} and \cite{KM01}, while four-dimensional
examples are studied in \cite{IJL05}.

Many of these homogeneous solutions exhibit pointed Gromov--Hausdorff
convergence to lower dimensional manifolds. For such solutions, it is not
possible to form a limit solution $(\mathcal{M}_{\infty}^{n},g_{\infty
}(t),x_{\infty})$ in a naive way. Nonetheless, an old conjecture of Hamilton
makes the following prediction.

\begin{conjecture}
[Hamilton]\label{Old}A Type-III solution should, in a suitable sense, behave
asymptotically like a (locally) homogeneous model.
\end{conjecture}

Evidence in favor of Hamilton's conjecture is found in prior work of the
authors, specifically \cite{HI93} and \cite{Knopf00}. As noted above, further
evidence comes from recent work of Lott \cite{Lott06}. Let $(\mathcal{M}%
^{n},g(t))$ be a solution of Ricci flow, and let $(\mathcal{M}^{n}%
,g_{s}(t),x)$ denote the rescaled pointed solution of Ricci flow defined for
$t\in\lbrack0,\infty)$ by $g_{s}(t)=\frac{1}{s}g(st)$.

\begin{theorem}
[Lott]\label{Lott1}Let $(\mathcal{M}^{n},g(t))$ be a Type-III solution of
Ricci flow. If the limit%
\[
(\mathcal{M}_{\infty}^{n},g_{\infty}(t),x_{\infty})=\lim_{s\rightarrow\infty
}(\mathcal{M}^{n},g_{s}(t),x)
\]
exists, then $(\mathcal{M}_{\infty}^{n},g_{\infty}(t))$ is an expanding Ricci soliton.
\end{theorem}

Bear in mind that this convergence is modulo (time-independent)
diffeomorphisms $\psi_{s}$, which one should not expect to converge as
$s\rightarrow\infty$. Note too that the word `expanding' refers to the sign of
$\alpha$ in (\ref{RSSE}), hence is compatible with collapse.

\smallskip

In dimension $n=3$, one can make a stronger statement, one that does not need
convergence as a hypothesis.

\begin{theorem}
[Lott]\label{Lott2}Let $(\mathcal{M}^{3},g(t))$ be a solution of Ricci flow on
a simply-connected homogeneous space $\mathcal{M}^{3}=G/K$, where $G$ is a
unimodular Lie group and $K$ is a compact isotropy subgroup. Then there exists
a limit%
\[
(\mathcal{M}_{\infty}^{3},g_{\infty}(t),x_{\infty})=\lim_{s\rightarrow\infty
}(\mathcal{M}^{3},g_{s}(t),x)
\]
which is an expanding homogeneous soliton (on a possibly different Lie group).
\end{theorem}

In light of Theorems~\ref{Lott1} and \ref{Lott2}, it is natural to refine
Conjecture~\ref{Old} by comparing the behavior of an arbitrary Type-III
solution to an expanding homogeneous soliton.

\subsection{Expanding homogeneous solitons}

Cheeger--Fukaya--Gromov show that any Riemannian manifold $(\mathcal{M}%
^{n},g)$ which collapses with bounded curvature is close to one with a
nilpotent structure. More precisely, there is for any $\varepsilon>0$ a metric
$g_{\varepsilon}$ on $\mathcal{M}^{n}$ which is $C^{1}$-close to $g$ such that
$g_{\varepsilon}$ admits a sheaf of nilpotent Lie algebras of local Killing
vector fields acting in the collapsed directions \cite{CFG92}. (Also see
\cite{CG86, CG90} and \cite{Fukaya87, Fukaya89}.) For this reason, when
looking for homogeneous models of Type-III behavior, it is natural to begin by
seeking expanding Ricci soliton structures on nilpotent Lie groups, which we
call \emph{Ricci nilsolitons}.

This search has an interesting history. Lauret \cite{Lauret01} approaches it
by asking a different question: is there a `best' $G$-invariant Riemannian
metric on a given homogeneous space $G$?\ An Einstein metric would be an
obvious candidate. Such metrics exist on some solvable Lie groups. In fact,
every known example of a noncompact, nonflat homogeneous space that admits an
Einstein metric is isomorphic to a solvable Lie group $S$ with a metric of
`standard type'. See \cite{Heber98} and \cite{Schueth2004}. (A left-invariant
metric $g$ on a solvable Lie group $S$, regarded as an inner product on the
Lie algebra $\mathfrak{s}$, is said to be of \emph{standard type} if the
orthogonal complement with respect to $g$ of the derived algebra
$[\mathfrak{s},\mathfrak{s}]$ forms an abelian subalgebra $\mathfrak{a}$ of
$\mathfrak{s}$.)

On the other hand, a noncommutative, nilpotent Lie group admits no Einstein
metric whatsoever. Indeed, by a result of Milnor, if the Lie algebra
$\mathfrak{g}$ of $G$ is nilpotent but noncommutative, then the Ricci
curvature of any left-invariant metric on $G$ has mixed sign \cite[Theorem
2.4]{Milnor76}. Since (\ref{RSSE}) reduces to $2\operatorname*{Rc}=-\alpha g$
when $X=0$, a Ricci soliton may be regarded as a generalization of an Einstein
metric, and thus provides a substitute candidate for the role of `best'
metric. Lauret proves the following result.

\begin{theorem}
[Lauret]\label{Lauret1}Let $N^{n}$ be a simply-connected nilpotent Lie group
with Lie algebra $\mathfrak{n}$ and a left-invariant metric $g$. Then
$(N^{n},g)$ admits a Ricci soliton structure if and only if $(\mathfrak{n},g)$
admits a standard metric solvable extension $(\mathfrak{s}=\mathfrak{a}%
\oplus\mathfrak{n},\tilde{g})$ such that the simply-connected solvable Lie
group $(S,\tilde{g})$ is Einstein.
\end{theorem}

If $\mathfrak{n}$ is a Lie algebra with inner product $g$, a \emph{standard
metric solvable extension} of $(\mathfrak{n},g)$ is $(\mathfrak{s},\tilde{g}%
)$, where $\mathfrak{s}=\mathfrak{a}\oplus\mathfrak{n}$ is a solvable Lie
algebra such that $[\cdot,\cdot]_{\mathfrak{s}}|_{\mathfrak{n}\times
\mathfrak{n}}=[\cdot,\cdot]_{\mathfrak{n}}$, and $\tilde{g}$ is an inner
product of standard type such that $\tilde{g}|_{\mathfrak{n}\times
\mathfrak{n}}=g$ and $[\mathfrak{s},\mathfrak{s}]_{\mathfrak{s}}%
=\mathfrak{n=}\mathfrak{a}^{\perp}$.

\medskip

As we noted above, the space of left-invariant metrics on $N^{n}$ has
dimension $n(n+1)/2$. However, a Ricci soliton structure, if one exists, is
essentially unique.

\begin{theorem}
[Lauret]\label{Lauret2}Let $N^{n}$ be a simply-connected nilpotent Lie group
with Lie algebra $\mathfrak{n}$. If $g$ and $g^{\prime}$ are Ricci soliton
metrics, then there exist $a>0$ and $\eta\in\operatorname*{Aut}(\mathfrak{n})$
such that $g^{\prime}=a\eta(g)$.
\end{theorem}

Although its proof is nonconstructive, Theorem~\ref{Lauret1} provides many
examples of existence or nonexistence. Any generalized Heisenberg group
\cite{BTV95} and many other two-step nilpotent Lie groups admit a Ricci
soliton structure. On the other hand, if $\mathfrak{n}$ is characteristically
nilpotent, then $N^{n}$ admits no such structure. (See \cite{Lauret01} and
\cite{LW06} for more examples.)

The first explicit constructions of Ricci soliton structures on nilpotent (or
more generally, solvable) Lie groups are obtained by Baird and Danielo
\cite{BD06} in dimension $n=3$ and independently by Lott \cite{Lott06} in
dimensions $n=3,4$. Remarkably, these are the first known examples of
nongradient soliton structures. See equations (\ref{dX3}), (\ref{dXSol}), and
(\ref{dX4}), below.

\smallskip

If $(N,g,X,\alpha)$ is an expanding Ricci soliton structure on a Lie group
that admits a compact quotient $\mathcal{M}^{n}=N/K$, what can one say about
soliton structures on $\mathcal{M}^{n}$? Because the metric $g$ is
homogeneous, it is compatible with any compact quotient. But the soliton
structure will never be compatible with compactification, because any compact
steady or expanding soliton is Einstein. (See \cite{Ivey93} or
\cite{Hamilton95}.) Accordingly, we call such compact quotients
\emph{pseudosolitons} --- compact manifolds with no soliton structure which
acquire one when lifted to an infinite cover. Apart from the Baird--Danielo
and Lott examples, the only other explicit example we know is $S^{1}\times
S^{n}$ with a product metric. This is not a soliton for any $n\geq2$, but it
is a quotient of $(\mathbb{R}\times S^{n},ds^{2}+2(n-1)g_{\operatorname*{can}%
},\operatorname*{grad}(s^{2}/4),-1)$. Pseudosolitons open the door to a more
complete understanding of the analytic behavior of Type-III solutions, and in
particular of Conjecture~\ref{Old}. That is to say, if Hamilton's expectation
is true, then nilpotent (or more generally, solvable) pseudosolitons should
act in some sense as stable models of Type-III behavior.

\begin{conjecture}
\label{New}Let $(\mathcal{M}^{n},g(t))$ be a solution of Ricci flow such that
$g(0)$ is sufficiently close (in a suitable norm) to a (locally) homogeneous
metric on (a quotient of) a nilpotent Lie group. Then there exists a limit
solution $(\mathcal{M}_{\infty}^{n},g_{\infty}(t),x_{\infty})$ which is a
nilsoliton.\footnote{One expects the nilsoliton to be a solution of Ricci flow
on $\mathcal{M}_{\infty}^{n}=\mathcal{\tilde{M}}^{n}$ obtained as a limit of
large-time asymptotic blow-downs, modulo diffeomorphisms, as in \cite{Lott06}%
.}
\end{conjecture}

If $(\mathcal{M}^{n},g(t))$ is a solution on a compact manifold, properties of
$(\mathcal{M}_{\infty}^{n},g_{\infty}(t),x_{\infty})$ should provide
information about $(\mathcal{M}^{n},g(t))$ via equivariant convergence theory,
for example, Lott's Riemannian groupoid results \cite{Lott06} and potential extensions.

\section{A modified Ricci flow and its linearization\label{Linear}}

\subsection{The flow}

Suppose one is given an immortal solution $(\mathcal{M}^{n},\hat{g}(t))$ of
Ricci flow%
\[
\frac{\partial}{\partial t}\hat{g}=-2\operatorname*{Rc}(\hat{g}),\qquad
\qquad(\beta<t<\infty),
\]
a constant $\alpha>0$, and a time-independent vector field $X$ on
$\mathcal{M}^{n}$. Define a scaling factor
\[
\sigma(t)=\alpha(t-\beta)
\]
and a one-parameter family of vector fields%
\[
Y(x,t)=\frac{1}{\sigma(t)}X(x).
\]
Let $\eta_{t}$ denote the family of diffeomorphisms of $\mathcal{M}^{n}$
generated by $Y$, so that%
\begin{equation}
\frac{\partial}{\partial t}\eta_{t}(x)=Y(\eta_{t}(x),t)=\frac{1}{\sigma
(t)}X(\eta_{t}(x)). \label{Diffeo}%
\end{equation}
Assume that $\eta_{t}$ exists for $t>\beta$. Because we are not assuming that
$\mathcal{M}^{n}$ is compact, this assumption requires verification in our
applications. (See Propositions~\ref{Exist3}, \ref{ExistSol}, and \ref{Exist4}
below.) Let%
\[
\tau(t)=\frac{1}{\alpha}\log(t-\beta),
\]
noting that $d\tau/dt=1/\sigma(t)$. Define a one-parameter family of metrics
$g(\tau)$ on $\mathcal{M}^{n}$ by%
\[
g(\tau)=\frac{1}{\sigma(\beta+e^{\alpha\tau})}(\eta_{\beta+e^{\alpha\tau}%
}^{-1})^{\ast}[\hat{g}(\beta+e^{\alpha\tau})].
\]
Observe that%
\[
\hat{g}(t)=\sigma(t)\eta_{t}^{\ast}\left[  g(\tau(t))\right]
\]
and that%
\begin{align*}
\frac{\partial}{\partial t}\hat{g}(t)  &  =\sigma^{\prime}(t)\eta_{t}^{\ast
}\left[  g(\tau(t))\right]  +\sigma(t)\eta_{t}^{\ast}\left\{  \mathcal{L}%
_{Y(t)}g(\tau(t))+\frac{d\tau}{dt}\left[  \frac{\partial}{\partial\tau}%
g(\tau)\right]  \right\} \\
&  =\eta_{t}^{\ast}\left\{  \alpha g(\tau)+\mathcal{L}_{X}g(\tau
)+\frac{\partial}{\partial\tau}g(\tau)\right\}  .
\end{align*}
Recalling that $\hat{g}(t)$ is a solution of Ricci flow, one sees that%
\[
\eta_{t}^{\ast}\left\{  -2\operatorname*{Rc}(g)\right\}  =-2\operatorname*{Rc}%
(\hat{g})=\eta_{t}^{\ast}\left\{  \alpha g(\tau)+\mathcal{L}_{X}g(\tau
)+\frac{\partial}{\partial\tau}g(\tau)\right\}  ,
\]
hence that $g(\tau)$ evolves by the modified Ricci flow%
\begin{equation}
\frac{\partial}{\partial\tau}g=-2\operatorname*{Rc}(g)-\mathcal{L}_{X}g-\alpha
g,\qquad\qquad(-\infty<\tau<\infty). \label{MRF}%
\end{equation}
The following observation then follows immediately from (\ref{RSSE}).

\begin{proposition}
If $(\mathcal{M}^{n},g,X,\alpha)$ is an expanding Ricci soliton structure,
then $g$ is a fixed point of the modified Ricci flow (\ref{MRF}).
\end{proposition}

\begin{remark}
This construction is easily adapted to shrinking $(\alpha<0)$ or steady
$(\alpha=0)$ solitons.
\end{remark}

\subsection{The linearization}

The linearization of (\ref{MRF}) is entirely standard. We recall it here for
the convenience of the reader.

Let $h$ be an arbitrary smooth $(2,0)$-tensor field with compact support on
$(\mathcal{M}^{n},g)$. Let $\{\tilde{g}(s):-\varepsilon<s<\varepsilon\}$ be a
one-parameter family of metrics with $\tilde{g}(0)=g$ and
\[
\frac{\partial}{\partial s}\tilde{g}|_{s=0}=h.
\]
Define $H=\operatorname*{tr}\!_{g}h$. It is well known that the linearization
of $-2\operatorname*{Rc}(g)$ at $g$ is%
\[
\frac{\partial}{\partial s}\left(  -2R_{ij}\right)  |_{s=0}=\Delta_{\ell
}h_{ij}+\nabla_{i}\nabla_{j}H+\nabla_{i}(\delta h)_{j}+\nabla_{j}(\delta
h)_{i},
\]
where $\Delta_{\ell}$ denotes the Lichnerowicz Laplacian%
\begin{equation}
\Delta_{\ell}h_{ij}=\Delta h_{ij}+2R_{ipqj}h^{pq}-R_{i}^{k}h_{kj}-R_{j}%
^{k}h_{ik}. \label{Lich}%
\end{equation}
We impose ellipticity by the DeTurck trick \cite{DT83, DT03}. Let $\{\tilde
{X}(s):-\varepsilon<s<\varepsilon\}$ be a one-parameter family of smooth
compactly supported vector fields with $\tilde{X}(0)=X$ and\thinspace
\footnote{If $\theta$ is a $1$-form, we denote by $\theta^{^{\sharp}}$ the
vector field metrically dual to $\theta$.}
\[
\frac{\partial}{\partial s}\tilde{X}|_{s=0}=W:=(\frac{1}{2}dH+\delta
h)^{\sharp}.
\]
Observing that\thinspace\footnote{This choice of $W$ is essentially equivalent
to the Bianchi gauge, as adopted in the elliptic context by Biquard
\cite{Biquard00} and others.}
\[
\frac{\partial}{\partial s}(\mathcal{L}_{\tilde{X}}\tilde{g})|_{s=0}%
=\mathcal{L}_{W}g+\mathcal{L}_{X}h,
\]
where $(\mathcal{L}_{W}g)_{ij}=\nabla_{i}\nabla_{j}H+\nabla_{i}(\delta
h)_{j}+\nabla_{j}(\delta h)_{i}$, one reaches the following conclusion.

\begin{proposition}
The linearization of the modified Ricci flow (\ref{MRF}) at a Ricci soliton
$(\mathcal{M}^{n},g,X,\alpha)$ is given by%
\[
\frac{\partial}{\partial\tau}h=Lh,
\]
where $L$ is the elliptic operator%
\begin{equation}
L:h\mapsto\Delta_{\ell}h-\mathcal{L}_{X}h-\alpha h. \label{L}%
\end{equation}

\end{proposition}

\emph{A priori, }$L$ is defined only on smooth compactly supported tensor
fields. Below, we will specify a domain that is useful for our purposes.

\begin{remark}
\label{BetterL}It is convenient to write $L$ in the alternative form%
\begin{equation}
L:h\mapsto\Delta_{\ell}h-\nabla_{X}h-\alpha h-\Xi(h), \label{LX}%
\end{equation}
where in coordinates, $\Xi(h)_{ij}=\nabla_{i}X^{k}h_{kj}+\nabla_{j}X^{k}%
h_{ki}$.
\end{remark}

\begin{example}
The Gaussian soliton is $(\mathbb{R}^{n},g,X,\alpha)$, where $g$ is the
standard flat metric, $\alpha\in\mathbb{R}$ is arbitrary, and
\[
X(x)=\operatorname*{grad}(-\frac{\alpha}{4}|x|^{2}).
\]
The linearization (\ref{L}) is simply%
\[
L:h\mapsto\Delta h-\nabla_{X}h,
\]
which is self adjoint in the weighted space $L^{2}(\mathbb{R}^{n}%
;\,e^{\frac{\alpha}{4}|x|^{2}}d\mu)$, where $d\mu=d\mu(g)$.
\end{example}

Critically, it turns out that the linearization is self adjoint in a weighted
space if and only if there is a \emph{gradient }soliton structure, i.e.~if and
only if $X=\operatorname*{grad}\varphi$ for some potential function
$\varphi:\mathcal{M}^{n}\rightarrow\mathbb{R}$.

\begin{remark}
\label{SelfAdjoint}Let $L$ be the operator (\ref{L}), let $\varphi$ be a
smooth function, and let $\left(  \cdot,\cdot\right)  _{\varphi}$ denote the
weighted $L^{2}$ inner product $\left(  u,v\right)  _{\varphi}=\int
_{\mathcal{M}^{n}}\left\langle u,v\right\rangle \,e^{-\varphi}d\mu$. Then one
has%
\begin{align*}
\left(  Lu,v\right)  _{\varphi}-\left(  u,Lv\right)  _{\varphi}  &  =2\left(
u,\nabla_{(X-\nabla\varphi)}v\right)  _{\varphi}+2\left(  dX_{\flat
},uv\right)  _{\varphi}\\
&  -\int_{\mathcal{M}^{n}}\left\langle u,v\right\rangle \left\{  \Delta
\varphi+\nabla_{X}\varphi-|\nabla\varphi|^{2}+\delta X\right\}  \,e^{-\varphi
}d\mu,
\end{align*}
where $(dX_{\flat})_{ij}=\nabla_{i}X_{j}-\nabla_{j}X_{i}$ in coordinates.
Clearly, the right-hand side above vanishes if $X=\operatorname*{grad}\varphi
$. The converse is left as an exercise.
\end{remark}

\section{Spectral bounds for known examples\label{Spectra}}

\subsection{$\operatorname*{nil}^{3}$ geometry\label{n3}}

We now describe the $\operatorname*{nil}^{3}$ soliton $(\mathcal{M}%
^{3},g,X,3)$ constructed in \cite{BD06}. (Also see \cite[Section~3.3.3]%
{Lott06}.) Because the exponential map of any connected, simply-connected,
nilpotent Lie group is a diffeomorphism, $\mathcal{M}^{3}$ is diffeomorphic to
$\mathbb{R}^{3}$. In standard coordinates $(x_{1},x_{2},x_{3})$ on
$\mathbb{R}^{3}$, consider the frame field $F=(F_{1},F_{2},F_{3})$ given by
\[
F_{1}=2\frac{\partial}{\partial x_{1}},\;F_{2}=2(\frac{\partial}{\partial
x_{2}}-x_{1}\frac{\partial}{\partial x_{3}}),\;F_{3}=2\frac{\partial}{\partial
x_{3}}.
\]
It is easy to check that all brackets $[F_{i},F_{j}]$ vanish except
$[F_{1},F_{2}]=-2F_{3}$. The connection is represented by the matrix%
\[
(\nabla_{F_{i}}F_{j})=%
\begin{pmatrix}
0 & -F_{3} & F_{2}\\
F_{3} & 0 & -F_{1}\\
F_{2} & -F_{1} & 0
\end{pmatrix}
.
\]
With respect to the dual field%
\[
\varphi^{1}=\frac{1}{2}\,dx_{1},\;\varphi^{2}=\frac{1}{2}\,dx_{2}%
,\;\varphi^{3}=\frac{1}{2}(x_{1}\,dx_{2}+\,dx_{3}),
\]
one may identify $g=g_{ij}\,\varphi^{i}\otimes\,\varphi^{j}$ with the matrix%
\[
(g_{ij})=%
\begin{pmatrix}
4 & 0 & 0\\
0 & 4 & 0\\
0 & 0 & 4
\end{pmatrix}
.
\]
Recalling the standard formula for the Riemannian curvature,
\begin{align*}
\left\langle R(X,Y)Y,X\right\rangle  &  =\frac{1}{4}\left\vert
(\operatorname{ad}X)^{\ast}Y+(\operatorname{ad}Y)^{\ast}X\right\vert
^{2}-\left\langle (\operatorname{ad}X)^{\ast}X,(\operatorname{ad}Y)^{\ast
}Y\right\rangle \\
&  -\frac{3}{4}\left\vert [X,Y]\right\vert ^{2}-\frac{1}{2}\left\langle
\left[  [X,Y],Y\right]  ,X\right\rangle -\frac{1}{2}\left\langle \left[
[Y,X],X\right]  ,Y\right\rangle ,
\end{align*}
we compute the Ricci curvature tensor $\operatorname*{Rc}=R_{ij}\,\varphi
^{i}\otimes\,\varphi^{j}$ of $g$ and identify it with%
\[
(R_{ij})=%
\begin{pmatrix}
-2 & 0 & 0\\
0 & -2 & 0\\
0 & 0 & 2
\end{pmatrix}
.
\]

Define the vector field%
\[
X=-\frac{1}{2}x_{1}F_{1}-\frac{1}{2}x_{2}F_{2}-(\frac{1}{2}x_{1}x_{2}%
+x_{3})F_{3}.
\]
A calculation shows that $\nabla X=\nabla_{i}X^{j}\,\varphi^{i}\otimes\,F_{j}$
has components%
\begin{equation}
(\nabla_{i}X^{j})=%
\begin{pmatrix}
-1 & -(\frac{1}{2}x_{1}x_{2}+x_{3}) & -\frac{1}{2}x_{2}\\
\frac{1}{2}x_{1}x_{2}+x_{3} & -1 & \frac{1}{2}x_{1}\\
\frac{1}{2}x_{2} & -\frac{1}{2}x_{1} & -2
\end{pmatrix}
. \label{dX3}%
\end{equation}
It follows that $2\operatorname*{Rc}(g)+\mathcal{L}_{X}g+3g=0$, which verifies
the Ricci soliton structure. Let $\xi$ denote the $1$-form metrically dual to
$X$. Recalling the standard identities%
\begin{align*}
d\xi(V,W)  &  =V\left\langle X,W\right\rangle -W\left\langle X,V\right\rangle
-\left\langle X,[V,W]\right\rangle \\
&  =\left\langle \nabla_{V}X,W\right\rangle -\left\langle \nabla
_{W}X,V\right\rangle ,
\end{align*}
one observes that $\xi$ is not closed, hence that the soliton is nongradient.

The following observation confirms that we may study this soliton by the
methods of Section~\ref{Linear}.

\begin{proposition}
\label{Exist3}For any `big bang' time $\beta<0$, the diffeomorphisms $\eta
_{t}:\mathcal{M}^{3}\rightarrow\mathcal{M}^{3}$ defined by%
\[
\eta_{t}(x_{1},x_{2},x_{3})=\left(  (\frac{\beta}{\beta-t})^{1/3}%
x_{1},\;(\frac{\beta}{\beta-t})^{1/3}x_{2},\;(\frac{\beta}{\beta-t}%
)^{2/3}x_{3}\right)
\]
exist for all $t>\beta$ and satisfy (\ref{Diffeo}), with $\eta_{0}%
=\operatorname*{id}$.
\end{proposition}

\begin{proof}
The result follows readily once one checks that in standard coordinates,%
\begin{equation}
X=-x_{1}\frac{\partial}{\partial x_{1}}-x_{2}\frac{\partial}{\partial x_{2}%
}-2x_{3}\frac{\partial}{\partial x_{3}}. \label{X3}%
\end{equation}

\end{proof}

The relationship between the soliton structure $(\mathcal{M}^{3},g,X,3)$ and
the evolving family of metrics whose existence is guaranteed by
Theorem~\ref{Lott2} is elucidated by the following.

\begin{example}
Consider the time-dependent frame $E(t):=FA(t)$, where%
\[
A(t)=%
\begin{pmatrix}
0 & 0 & a(t)\\
0 & a(t) & 0\\
-2a^{2}(t) & 0 & 0
\end{pmatrix}
\qquad\text{and}\qquad a(t)=\sqrt{\frac{1}{12}t^{-2/3}}.
\]
With respect to $E(t)$, one obtains the identification%
\[
(g_{ij})_{E(t)}=%
\begin{pmatrix}
\frac{1}{9}t^{-4/3} & 0 & 0\\
0 & \frac{1}{3}t^{-2/3} & 0\\
0 & 0 & \frac{1}{3}t^{-2/3}%
\end{pmatrix}
.
\]
Observe that we may identify $3tg_{E(t)}$ with Lott's limit solution
\cite[Formula~(3.18)]{Lott06}
\[
g_{\infty}(t)=\frac{1}{3t^{1/3}}(\theta^{1}\otimes\,\theta^{1})+t^{1/3}%
(\theta^{2}\otimes\,\theta^{2})+t^{1/3}(\theta^{3}\otimes\,\theta^{3})
\]
given in coordinates $(x,y,z)$, with $\theta^{1}=dx+\frac{1}{2}y\,dz-\frac
{1}{2}z\,dy$, $\theta^{2}=dy$, and $\theta^{3}=dz$.
\end{example}

Now let $h=h_{ij}\,\varphi^{i}\otimes\,\varphi^{j}$ be a smooth compactly
supported tensor field. Let $H=\operatorname*{tr}\!_{g}h$, noting that%
\[
H=\frac{1}{4}(h_{11}+h_{22}+h_{33})
\]
and%
\[
|h|^{2}=\frac{1}{16}\left\{  h_{11}^{2}+h_{22}^{2}+h_{33}^{2}+2(h_{12}%
^{2}+h_{13}^{2}+h_{23}^{2})\right\}
\]
pointwise. Using Remark~\ref{BetterL} in conjunction with the alternative
formula%
\begin{equation}
\Delta_{\ell}h_{ij}=\Delta h_{ij}+(Rh_{ij}+2HR_{ij})-3(R_{i}^{k}h_{kj}%
+R_{j}^{k}h_{ik})+(2\left\langle \operatorname*{Rc},h\right\rangle -RH)g_{ij}
\label{Lich3}%
\end{equation}
for the Lichnerowicz Laplacian (\ref{Lich}), valid when $n=3$, one computes
that%
\begin{align*}
\int\left\langle Lh,h\right\rangle \,d\mu &  =\int\left\langle \Delta
h,h\right\rangle \,d\mu-\int\left\langle \nabla_{X}h,h\right\rangle \,d\mu\\
&  +\int(R-\alpha)\left\vert h\right\vert ^{2}\,d\mu+\int(4\left\langle
\operatorname*{Rc},h\right\rangle -RH)H\,d\mu\\
&  -6\int R_{i}^{j}h_{j}^{k}h_{k}^{i}\,d\mu-2\int\nabla_{i}X^{j}h_{j}^{k}%
h_{k}^{i}\,d\mu
\end{align*}
where all indices are with respect to the orthogonal frame field $F$.

Clearly, $\int\left\langle \Delta h,h\right\rangle \,d\mu=-\left\Vert \nabla
h\right\Vert ^{2}$. Observing that $\delta X=4$, one integrates by parts to
get
\[
-\int\left\langle \nabla_{X}h,h\right\rangle \,d\mu=-\frac{1}{2}\int(\delta
X)\left\langle h,h\right\rangle \,d\mu=-2\left\Vert h\right\Vert ^{2}.
\]
Because $R=-1/2$ and $\alpha=3$, one has $\int(R-\alpha)\left\vert
h\right\vert ^{2}\,d\mu=-\frac{7}{2}\left\Vert h\right\Vert ^{2}$. The
pointwise calculation $\left\langle \operatorname*{Rc},h\right\rangle
=\frac{1}{8}(-h_{11}-h_{22}+h_{33})$ shows that the fourth term reduces to
\[
\int(4\left\langle \operatorname*{Rc},h\right\rangle -RH)H\,d\mu=-\frac{3}%
{2}\left\Vert H\right\Vert ^{2}+\int h_{33}H\,d\mu.
\]
Calculating pointwise that $R_{i}^{j}h_{j}^{k}h_{k}^{i}=\frac{1}{32}%
(-h_{11}^{2}-h_{22}^{2}+h_{33}^{2}-2h_{12}^{2})$, one gets%
\[
-6\int R_{i}^{j}h_{j}^{k}h_{k}^{i}\,d\mu=\frac{3}{16}\int(h_{11}^{2}%
+h_{22}^{2}-h_{33}^{2}+2h_{12}^{2})\,d\mu.
\]
Computing $\nabla_{i}X^{j}h_{j}^{k}h_{k}^{i}=-|h|^{2}-\frac{1}{16}(h_{33}%
^{2}+h_{13}^{2}+h_{23}^{2})$ pointwise then gives%
\[
-2\int\nabla_{i}X^{j}h_{j}^{k}h_{k}^{i}\,d\mu=2\left\Vert h\right\Vert
^{2}+\frac{1}{8}\int(h_{33}^{2}+h_{13}^{2}+h_{23}^{2})\,d\mu.
\]
Putting these all together, one has%
\[
\int\left\langle Lh,h\right\rangle \,d\mu=-\left\Vert \nabla h\right\Vert
^{2}-\frac{1}{2}\left\Vert h\right\Vert ^{2}-\frac{3}{2}\left\Vert
H\right\Vert ^{2}+\int h_{33}H\,d\mu-\frac{1}{4}\int(h_{33}^{2}+h_{13}%
^{2}+h_{23}^{2})\,d\mu.
\]
Since $|h_{33}H|\leq\frac{1}{4}h_{33}^{2}+H^{2}$ by weighted Cauchy--Schwarz,
it follows that
\[
\int\left\langle Lh,h\right\rangle \,d\mu\leq-\left\Vert \nabla h\right\Vert
^{2}-\frac{1}{2}(\left\Vert h\right\Vert ^{2}+\left\Vert H\right\Vert ^{2}).
\]
By \cite{EF98}, smooth compactly supported tensor fields are dense in $L^{2}$.
Hence we have proved the following.

\begin{proposition}
\label{Stable3}The linearization of the modified Ricci flow (\ref{MRF}) at the
$\operatorname*{nil}^{3}$ soliton $(\mathcal{M}^{3},g,X,3)$ is strictly
linearly stable and satisfies%
\[
\int\left\langle Lh,h\right\rangle \,d\mu\leq-\left\Vert \nabla h\right\Vert
^{2}-\omega\left\Vert h\right\Vert ^{2}%
\]
with $\omega=1/2$.
\end{proposition}

\subsection{$\operatorname*{sol}^{3}$ geometry}

A $\operatorname*{sol}^{3}$ soliton $(\mathcal{M}^{3},g,X,4)$ has been found
by Baird--Danielo \cite{BD06} and independently by Lott \cite[Section~3.3.2]%
{Lott06}. $\mathcal{M}^{3}$ is diffeomorphic to $\mathbb{R}^{3}$. In standard
coordinates $(x_{1},x_{2},x_{3})$, define a frame field $F=(F_{1},F_{2}%
,F_{3})$ by%
\[
F_{1}=2\frac{\partial}{\partial x_{1}},\;F_{2}=2(e^{-x_{1}}\frac{\partial
}{\partial x_{2}}+e^{x_{1}}\frac{\partial}{\partial x_{3}}),\;F_{3}%
=2(e^{-x_{1}}\frac{\partial}{\partial x_{2}}-e^{x_{1}}\frac{\partial}{\partial
x_{3}}).
\]
The bracket relations are $[F_{1},F_{2}]=-2F_{3}$, $[F_{2},F_{3}]=0$, and
$[F_{3},F_{1}]=2F_{2}$. The connection is represented by%
\[
(\nabla_{F_{i}}F_{j})=%
\begin{pmatrix}
0 & 0 & 0\\
2F_{3} & 0 & -4F_{1}\\
2F_{2} & -4F_{1} & 0
\end{pmatrix}
.
\]
With respect to the dual field%
\[
\varphi^{1}=\frac{1}{2}\,dx_{1},\;\varphi^{2}=\frac{1}{4}(e^{x_{1}}%
\,dx_{2}+e^{-x_{1}}\,dx_{3}),\;\varphi^{3}=\frac{1}{4}(e^{x_{1}}%
\,dx_{2}-e^{-x_{1}}\,dx_{3}),
\]
we identify the metric $g=g_{ij}\,\varphi^{i}\otimes\,\varphi^{j}$ with the
matrix%
\[
(g_{ij})=%
\begin{pmatrix}
4 & 0 & 0\\
0 & 8 & 0\\
0 & 0 & 8
\end{pmatrix}
,
\]
and its Ricci tensor $\operatorname*{Rc}=R_{ij}\,\varphi^{i}\otimes
\,\varphi^{j}$ with%
\[
(R_{ij})=%
\begin{pmatrix}
-8 & 0 & 0\\
0 & 0 & 0\\
0 & 0 & 0
\end{pmatrix}
.
\]

Given any $\gamma\in\mathbb{R}$, define a vector field%
\[
X=\gamma\left(  -F_{1}-e^{-x_{1}}x_{3}F_{2}+e^{-x_{1}}x_{3}F_{3}\right)
+(1-\gamma)\left(  F_{1}-e^{x_{1}}x_{2}F_{2}-e^{x_{1}}x_{2}F_{3}\right)  .
\]
A calculation shows that the components of $\nabla X=\nabla_{i}X^{j}%
\,\varphi^{i}\otimes\,F_{j}$ are%
\begin{equation}
(\nabla_{i}X^{j})=\gamma%
\begin{pmatrix}
0 & 2e^{-x_{1}}x_{3} & -2e^{-x_{1}}x_{3}\\
-4e^{-x_{1}}x_{3} & -2 & 0\\
4e^{-x_{1}}x_{3} & 0 & -2
\end{pmatrix}
+(1-\gamma)%
\begin{pmatrix}
0 & -2e^{x_{1}}x_{2} & -2e^{x_{1}}x_{2}\\
4e^{x_{1}}x_{2} & -2 & 0\\
4e^{x_{1}}x_{2} & 0 & -2
\end{pmatrix}
. \label{dXSol}%
\end{equation}
It follows that $2\operatorname*{Rc}(g)+\mathcal{L}_{X}g+4g=0$, which verifies
the Ricci soliton structure. As in Section~\ref{n3}, one observes that the
soliton structure is nongradient.

Notice that in standard coordinates,%
\begin{equation}
X=\gamma\left(  -2\frac{\partial}{\partial x_{1}}-4x_{3}\frac{\partial
}{\partial x_{3}}\right)  +(1-\gamma)\left(  2\frac{\partial}{\partial x_{1}%
}-4x_{2}\frac{\partial}{\partial x_{2}}\right)  . \label{XSol}%
\end{equation}
We verify the existence of suitable diffeomorphisms for the case $\gamma=1/2$,
leaving the general case for an interested reader.

\begin{proposition}
\label{ExistSol}Take $\gamma=1/2$ in (\ref{XSol}). Then for any `big bang'
time $\beta<0$, the diffeomorphisms $\eta_{t}:\mathcal{M}^{3}\rightarrow
\mathcal{M}^{3}$ defined by%
\[
\eta_{t}(x_{1},x_{2},x_{3})=\left(  x_{1},\;(\frac{\beta}{\beta-t})^{1/2}%
x_{2},\;(\frac{\beta}{\beta-t})^{1/2}x_{3}\right)
\]
exist for all $t>\beta$ and satisfy (\ref{Diffeo}), with $\eta_{0}%
=\operatorname*{id}$.
\end{proposition}

The relationship between the soliton structure $(\mathcal{M}^{3},g,X,4)$ and
the evolving family of metrics given by Theorem~\ref{Lott2} can be seen as follows.

\begin{example}
Consider the time-dependent frame $E(t)=FA(t)$ given by%
\[
A(t)=%
\begin{pmatrix}
0 & -\frac{1}{2} & 0\\
a(t) & 0 & 0\\
0 & 0 & a(t)
\end{pmatrix}
\qquad\text{and}\qquad a(t)=\sqrt{\frac{1}{32}t}.
\]
With respect to the frame $E(t)$, one has the identification%
\[
(g_{ij})_{E(t)}=%
\begin{pmatrix}
\frac{1}{4}t^{-1} & 0 & 0\\
0 & 1 & 0\\
0 & 0 & \frac{1}{4}t^{-1}%
\end{pmatrix}
.
\]
Observe that we may identify $4tg_{E(t)}$ with Lott's limit solution
\cite[Formula~(3.9)]{Lott06}%
\[
g_{\infty}(t)=(\theta^{1}\otimes\,\theta^{1})+4t(\theta^{2}\otimes\,\theta
^{2})+(\theta^{3}\otimes\,\theta^{3})
\]
given in coordinates $(x,y,z)$, with $\theta^{1}+\theta^{3}=e^{-z}\,dx$,
$\theta^{1}-\theta^{3}=e^{z}\,dy$, and $\theta^{2}=dz$.
\end{example}

Now let $h=h_{ij}\,\varphi^{i}\otimes\,\varphi^{j}$ be a smooth compactly
supported tensor field. Let $H=\operatorname*{tr}\!_{g}h$, noting that%
\[
H=\frac{1}{8}(2h_{11}+h_{22}+h_{33})
\]
and%
\[
|h|^{2}=\frac{1}{64}\left\{  4h_{11}^{2}+h_{22}^{2}+h_{33}^{2}+4h_{12}%
^{2}+4h_{13}^{2}+2h_{23}^{2}\right\}
\]
pointwise. As in Section~\ref{n3}, we proceed to evaluate%
\begin{align*}
\int\left\langle Lh,h\right\rangle \,d\mu &  =\int\left\langle \Delta
h,h\right\rangle \,d\mu-\int\left\langle \nabla_{X}h,h\right\rangle \,d\mu\\
&  +\int(R-\alpha)\left\vert h\right\vert ^{2}\,d\mu+\int(4\left\langle
\operatorname*{Rc},h\right\rangle -RH)H\,d\mu\\
&  -6\int\left\langle \operatorname*{Rc},h^{2}\right\rangle \,d\mu-2\int
\nabla_{i}X^{j}h_{j}^{k}h_{k}^{i}\,d\mu.
\end{align*}
Since $\alpha=4=\delta X$, $R=-2$, $\left\langle \operatorname*{Rc}%
,h\right\rangle =-\frac{1}{2}h_{11}$, and $\left\langle \operatorname*{Rc}%
,h^{2}\right\rangle =-\frac{1}{16}(2h_{11}^{2}+h_{12}^{2}+h_{13}^{2})$, one
finds that%
\begin{equation}
\int\left\langle Lh,h\right\rangle \,d\mu=-\left\Vert \nabla h\right\Vert
^{2}-\frac{1}{32}\int P\,d\mu, \label{Sol1}%
\end{equation}
where%
\[
P:=32Hh_{11}-4h_{11}^{2}+(h_{22}-h_{33})^{2}+4h_{23}^{2}.
\]
At a point where $h=\operatorname*{diag}(a,-2a,-2a)$, one has $\frac{1}%
{32}P=-\frac{3}{8}a^{2}=-2|h|^{2}$, so we must work harder.

As in \cite{Koiso79}, one defines a $(3,0)$-tensor $T=T(h)$ by%
\[
T_{ijk}:=\nabla_{k}h_{ij}-\nabla_{i}h_{jk}%
\]
and computes that%
\[
\left\Vert \nabla h\right\Vert ^{2}=\left\Vert \delta h\right\Vert ^{2}%
+\frac{1}{2}\left\Vert T\right\Vert ^{2}+\int\left\{  R_{ijk\ell}h^{i\ell
}h^{jk}-R_{i}^{k}h_{jk}h^{ij}\right\}  \,d\mu.
\]
The decomposition of the Riemann curvature tensor in dimension three,%
\[
R_{ijk\ell}=R_{i\ell}g_{jk}+R_{jk}g_{i\ell}-R_{ik}g_{j\ell}-R_{j\ell}%
g_{ik}-\frac{1}{2}R(g_{i\ell}g_{jk}-g_{ik}g_{j\ell}),
\]
implies the pointwise identity%
\[
R_{ijk\ell}h^{i\ell}h^{jk}-R_{i}^{k}h_{jk}h^{ij}=2\left\langle
\operatorname*{Rc},h\right\rangle H-3\left\langle \operatorname*{Rc}%
,h^{2}\right\rangle -\frac{1}{2}RH^{2}+\frac{1}{2}R|h|^{2}.
\]
Using this and discarding the $\left\Vert T\right\Vert ^{2}$ term, one can
improve (\ref{Sol1}), obtaining%
\[
\int\left\langle Lh,h\right\rangle \,d\mu\leq-\left\Vert \delta h\right\Vert
^{2}-\frac{1}{64}\int Q\,d\mu,
\]
where%
\[
Q:=16h_{11}^{2}+h_{22}^{2}+h_{33}^{2}+8h_{12}^{2}+8h_{13}^{2}+6h_{23}%
^{2}+(h_{22}-h_{33})^{2}+4h_{11}(h_{22}+h_{33}).
\]
Recall that $64|h|^{2}=4h_{11}^{2}+h_{22}^{2}+h_{33}^{2}+4h_{12}^{2}%
+4h_{13}^{2}+2h_{23}^{2}$. Recall too that for any $\varepsilon>0$, one has
$|4h_{11}(h_{22}+h_{33})|\leq4\varepsilon h_{11}^{2}+\frac{2}{\varepsilon
}(h_{22}^{2}+h_{33}^{2})$. The choice $\varepsilon=(3+\sqrt{17})/2$ solves
$(16-4\varepsilon)/4=1-2/\varepsilon$, hence is optimal. Thus we obtain the
following result.

\begin{proposition}
\label{StableSol}The linearization of the modified Ricci flow (\ref{MRF}) at
the $\operatorname*{sol}^{3}$ soliton $(\mathcal{M}^{3},g,X,4)$ is strictly
linearly stable and satisfies%
\[
\int\left\langle Lh,h\right\rangle \,d\mu\leq-\left\Vert \delta h\right\Vert
^{2}-\omega\left\Vert h\right\Vert ^{2}%
\]
with $\omega=(5-\sqrt{17})/2>0$.
\end{proposition}

\subsection{$\operatorname*{nil}^{4}$ geometry}

We now describe the soliton $(\mathcal{M}^{4},g,X,3)$ introduced in
\cite[Section~3.4.9]{Lott06}. (Compare \cite[Section~A.6]{IJL05}.) In standard
coordinates $(x_{1},x_{2},x_{3},x_{4})$ on $\mathcal{M}^{4}\approx
\mathbb{R}^{4}$, consider the frame field%
\[
F_{1}=\frac{\partial}{\partial x_{1}},\;F_{2}=\frac{\partial}{\partial x_{2}%
},\;F_{3}=\frac{\partial}{\partial x_{3}},\;F_{4}=x_{1}\frac{\partial
}{\partial x_{2}}+x_{2}\frac{\partial}{\partial x_{3}}+\frac{\partial
}{\partial x_{4}}.
\]
It is easy to check that $[F_{1},F_{4}]=F_{2}$ and $[F_{2},F_{4}]=F_{3}$ with
all other brackets vanishing. The dual basis is%
\[
\varphi^{1}=dx_{1},\;\varphi^{2}=dx_{2}-x_{1}dx_{4},\;\varphi^{3}=dx_{3}%
-x_{2}dx_{4},\;\varphi^{4}=dx_{4},
\]
and the connection is specified by%
\[
(\nabla_{F_{i}}F_{j})=\frac{1}{2}%
\begin{pmatrix}
0 & -F_{4} & 0 & F_{2}\\
-F_{4} & 0 & -F_{4} & F_{1}+F_{3}\\
0 & -F_{4} & 0 & F_{2}\\
-F_{2} & F_{1}-F_{3} & F_{2} & 0
\end{pmatrix}
.
\]
Let $g=g_{ij}\,\varphi^{i}\otimes\,\varphi^{j}$ denote Lott's metric
$g_{\infty}(\frac{1}{3})$, where by \cite[Formula~(3.66)]{Lott06},%
\[
g_{\infty}(t)=3^{1/3}t^{1/3}(\varphi^{1}\otimes\,\varphi^{1})+(\varphi
^{2}\otimes\,\varphi^{2})+3^{-1/3}t^{-1/3}(\varphi^{3}\otimes\,\varphi
^{3})+3^{2/3}t^{2/4}(\varphi^{4}\otimes\,\varphi^{4}).
\]
Then $g$ corresponds to the matrix%
\[
g=%
\begin{pmatrix}
1 & 0 & 0 & 0\\
0 & 1 & 0 & 0\\
0 & 0 & 1 & 0\\
0 & 0 & 0 & 1
\end{pmatrix}
,
\]
and its Ricci curvature tensor $\operatorname*{Rc}=R_{ij}\,\varphi^{i}%
\otimes\,\varphi^{j}$ corresponds to%
\[
\operatorname*{Rc}=%
\begin{pmatrix}
-1/2 & 0 & 0 & 0\\
0 & 0 & 0 & 0\\
0 & 0 & 1/2 & 0\\
0 & 0 & 0 & 1
\end{pmatrix}
.
\]

Consider the vector field%
\[
X=-2x_{1}F_{1}+(-3x_{2}+x_{1}x_{4})F_{2}+(-4x_{3}+x_{2}x_{4})F_{3}-x_{4}%
F_{4}.
\]
Calculating%
\begin{equation}
(\nabla_{i}X^{j})=%
\begin{pmatrix}
-2 & \frac{1}{2}x_{4} & 0 & \frac{1}{2}(3x_{2}-x_{1}x_{4})\\
-\frac{1}{2}x_{4} & -3 & \frac{1}{2}x_{4} & x_{1}+2x_{3}-\frac{1}{2}x_{2}%
x_{4}\\
0 & -\frac{1}{2}x_{4} & -4 & \frac{1}{2}(3x_{2}-x_{1}x_{4})\\
\frac{1}{2}(x_{1}x_{4}-3x_{2}) & \frac{1}{2}x_{2}x_{4}-x_{1}-2x_{3} & \frac
{1}{2}(x_{1}x_{4}-3x_{2}) & -1
\end{pmatrix}
, \label{dX4}%
\end{equation}
we have verified the nongradient Ricci soliton structure $2\operatorname*{Rc}%
(g)+\mathcal{L}_{X}g+3g=0$.

\begin{proposition}
\label{Exist4}For any `big bang' time $\beta<0$, the diffeomorphisms $\eta
_{t}:\mathcal{M}^{4}\rightarrow\mathcal{M}^{4}$ defined by%
\[
\eta_{t}(x_{1},x_{2},x_{3},x_{4})=\left(  (\frac{\beta}{\beta-t})^{2/3}%
x_{1},\;(\frac{\beta}{\beta-t})x_{2},\;(\frac{\beta}{\beta-t})^{4/3}%
x_{3},\;(\frac{\beta}{\beta-t})^{1/3}x_{4}\right)
\]
exist for all $t>\beta$ and satisfy (\ref{Diffeo}), with $\eta_{0}%
=\operatorname*{id}$.
\end{proposition}

\begin{proof}
The result follows immediately from the observation that in standard
coordinates,%
\begin{equation}
X=-2x_{1}\frac{\partial}{\partial x_{1}}-3x_{2}\frac{\partial}{\partial x_{2}%
}-4x_{3}\frac{\partial}{\partial x_{3}}-x_{4}\frac{\partial}{\partial x_{4}}.
\label{X4}%
\end{equation}

\end{proof}

Let $h=h_{ij}\,\varphi^{i}\otimes\,\varphi^{j}$ be a smooth compactly
supported tensor field. By formulas (\ref{LX}) and (\ref{Lich}), one has%
\begin{align*}
\int\left\langle Lh,h\right\rangle \,d\mu &  =\int\left\langle \Delta
h,h\right\rangle \,d\mu-\int\left\langle \nabla_{X}h,h\right\rangle \,d\mu\\
&  +2\int R_{ijk\ell}h^{i\ell}h^{jk}\,d\mu-2\int R_{i}^{j}h_{j}^{k}h_{k}%
^{i}\,d\mu\\
&  -3\left\Vert h\right\Vert ^{2}-2\int\nabla_{i}X^{j}h_{j}^{k}h_{k}^{i}%
\,d\mu,
\end{align*}
where all indices are with respect to the orthonormal frame $F$. Since formula
(\ref{Lich3}) is not available, one must calculate the full curvature tensor
to proceed. So define $R_{ijk}:=R(F_{i},F_{j})F_{k}=\nabla_{F_{i}}%
(\nabla_{F_{j}}F_{k})-\nabla_{F_{j}}(\nabla_{F_{i}}F_{k})-\nabla_{\lbrack
F_{i},F_{j}]}F_{k}$. Then it is straightforward to compute that all
nonvanishing curvature components are determined by%
\begin{align*}
\frac{1}{4}F_{2}  &  =R_{211}=R_{213}=R_{231}=R_{233}=\frac{1}{2}R_{424}\\
\frac{1}{4}F_{4}  &  =\frac{1}{3}R_{141}=R_{143}=\frac{1}{2}R_{242}%
=R_{341}=R_{433}\\
\frac{1}{4}F_{1}+\frac{1}{4}F_{3}  &  =R_{122}=R_{322}\\
\frac{1}{4}F_{1}-\frac{1}{4}F_{3}  &  =R_{434}\\
\frac{3}{4}F_{1}+\frac{1}{4}F_{3}  &  =R_{414}.
\end{align*}
Further pointwise calculations show that%
\begin{align*}
2R_{ijk\ell}h^{i\ell}h^{jk}  &  =-h_{12}^{2}+3h_{14}^{2}-h_{23}^{2}%
+2h_{24}^{2}-h_{34}^{2}\\
&  +h_{11}h_{22}-3h_{11}h_{44}-2h_{12}h_{23}+2h_{13}h_{22}-2h_{13}h_{44}\\
&  +2h_{14}h_{34}+h_{22}h_{33}-2h_{22}h_{44}+h_{33}h_{44}%
\end{align*}
and%
\[
2R_{i}^{j}h_{j}^{k}h_{k}^{i}=-h_{11}^{2}+h_{33}^{2}-2h_{44}^{2}-h_{12}%
^{2}-3h_{14}^{2}+h_{23}^{2}-2h_{24}^{2}-h_{34}^{2}%
\]
and
\begin{align*}
\nabla_{i}X^{j}h_{j}^{k}h_{k}^{i}  &  =-2h_{11}^{2}-3h_{22}^{2}-4h_{33}%
^{2}-h_{44}^{2}\\
&  -5h_{12}^{2}-6h_{13}^{2}-3h_{14}^{2}-4h_{24}^{2}-7h_{23}^{2}-5h_{34}^{2}.
\end{align*}
Because $\delta X=10$, one has $-\int\left\langle \nabla_{X}h,h\right\rangle
\,d\mu=-5\left\Vert h\right\Vert ^{2}$. Collecting terms, one then concludes
that%
\[
\int\left\langle Lh,h\right\rangle \,d\mu=-\left\Vert \nabla h\right\Vert
^{2}+\int Q\,d\mu,
\]
where%
\begin{align*}
Q  &  :=-3h_{11}^{2}-2h_{22}^{2}-h_{33}^{2}-4h_{44}^{2}-6h_{12}^{2}%
-4h_{13}^{2}-4h_{14}^{2}-4h_{23}^{2}-4h_{24}^{2}-6h_{34}^{2}\\
&  +h_{11}h_{22}-3h_{11}h_{44}-2h_{12}h_{23}+2h_{13}h_{22}-2h_{13}h_{44}\\
&  +2h_{14}h_{34}+h_{22}h_{33}-2h_{22}h_{44}+h_{33}h_{44}.
\end{align*}

\begin{proposition}
\label{Stable4}The linearization of the modified Ricci flow (\ref{MRF}) at the
$\operatorname*{nil}^{4}$ soliton $(\mathcal{M}^{4},g,X,3)$ is strictly
linearly stable and satisfies%
\[
\int\left\langle Lh,h\right\rangle \,d\mu\leq-\left\Vert \nabla h\right\Vert
^{2}-\omega\left\Vert h\right\Vert ^{2}%
\]
for some $\omega>0.0057$.
\end{proposition}

\begin{proof}
For constants $A,B,C,D,E,F>0$ to be chosen, weighted Cauchy--Schwarz implies
that $|h_{11}h_{22}|\leq\frac{1}{4A}h_{11}^{2}+Ah_{22}^{2}$, $2|h_{13}%
h_{22}|\leq\frac{1}{2B}h_{13}^{2}+2Bh_{22}^{2}$, $2|h_{13}h_{44}|\leq\frac
{1}{2C}h_{13}^{2}+2Ch_{44}^{2}$, $|h_{22}h_{33}|\leq Dh_{22}^{2}+\frac{1}%
{4D}h_{33}^{2}$, $2|h_{22}h_{44}|\leq2Eh_{22}^{2}+\frac{1}{2E}h_{44}^{2}$, and
$|h_{33}h_{44}|\leq Fh_{33}^{2}+\frac{1}{4F}h_{44}^{2}$. This leads to the
estimate%
\begin{align*}
Q  &  \leq(-\frac{3}{2}+\frac{1}{4A})h_{11}^{2}+(-2+A+2B+D+2E)h_{22}^{2}\\
&  +(-1+\frac{1}{4D}+F)h_{33}^{2}+(-\frac{5}{2}+2C+\frac{1}{2E}+\frac{1}%
{4F})h_{44}^{2}\\
&  -5h_{12}^{2}+(-4+\frac{1}{2B}+\frac{1}{2C})h_{13}^{2}-3h_{14}^{2}%
-3h_{23}^{2}-4h_{24}^{2}-5h_{34}^{2}.
\end{align*}
The choices $A=0.17$, $B=0.27$, $C=0.24$, $D=0.35$, $E=0.46$, and $F=0.28$
\emph{(which are not optimal!) }yield the estimate $Q\leq-\omega|h|^{2}$.
\end{proof}

\section{Generation of strongly continuous semigroups\label{Semigroups}}

Our approach in this section incorporates ideas of Rabier \cite{Rabier05},
adapting them to tensor-valued operators on noncompact manifolds. Roughly
speaking, the idea is to approximate an unbounded operator by a sequence of
bounded operators, obtain uniform estimates for this sequence, and then use
the fact that all but finitely many of the approximants agree with the
original operator when paired against any tensor field with compact support.

Let $(\mathcal{M}^{n},g,X,\alpha)$ be the $\operatorname*{nil}^{3}$,
$\operatorname*{sol}^{3}$, or $\operatorname*{nil}^{4}$ soliton discussed
above. For each $i,j,k\in\mathbb{N}$, we denote by $W_{i,j}^{k}\equiv
W_{i,j}^{k,2}$ the Sobolev space of complex-valued $(i,j)$-tensor fields on
$\mathcal{M}^{n}$, obtained as the completion of $C_{0}^{\infty}(T_{i}%
^{j}\mathcal{M}^{n};\mathbb{C})$ with respect to the norm
\[
\left\Vert u\right\Vert _{k}^{2}=\sum_{\ell=0}^{k}\int_{\mathcal{M}^{n}%
}\left\vert \nabla^{\ell}u\right\vert ^{2}\,d\mu=\sum_{\ell=0}^{k}\left(
\nabla^{\ell}u,\nabla^{\ell}u\right)  .
\]
To avoid notational prolixity, we write each $L_{i,j}^{2}=W_{i,j}^{0}$ inner
product as $\left(  \cdot,\cdot\right)  $; the relevant values of $i,j$ should
be clear from the context. We denote the real and imaginary parts of a complex
number $z$ by $\operatorname*{Re}(z)$ and $\operatorname{Im}(z)$, respectively.

\begin{lemma}
\label{LM}Let $\xi\in W_{0,1}^{1}$ be a real vector field with $\left\Vert
\delta\xi\right\Vert _{\infty}<\infty$, and let $\zeta$ be a real
$(2,2)$-tensor field with $\left\Vert \zeta\right\Vert _{\infty}<\infty$.
Assume there exists $\omega\in\mathbb{R}$ such that for all $v\in W_{2,0}^{1}%
$, one has
\begin{equation}
\left(  \frac{1}{2}(\delta\xi)v-\zeta v,\bar{v}\right)  \geq\omega\left\Vert
v\right\Vert _{0}^{2}. \label{Axiom}%
\end{equation}
Then for each $\lambda>-\omega$ and each $f\in L_{2,0}^{2}$, there exists a
unique $u\in W_{2,0}^{1}$ with norm%
\[
\left\Vert u\right\Vert _{1}^{2}\leq\frac{1}{\lambda+\omega}\left\Vert
f\right\Vert _{0}^{2}%
\]
such that $u$ solves $\beta_{\lambda}(u,v)=\left(  f,\bar{v}\right)  $ for all
$v\in W_{2,0}^{1}$, where%
\begin{equation}
\beta_{\lambda}(u,v):=\lambda\left(  u,\bar{v}\right)  +\left(  \nabla
u,\nabla\bar{v}\right)  +\left(  \nabla_{\xi}u,\bar{v}\right)  -\left(  \zeta
u,\bar{v}\right)  .
\end{equation}

\end{lemma}

\begin{proof}
For each $\lambda\in\mathbb{R}$, $\beta_{\lambda}$ is a continuous
sesquilinear form on $W_{2,0}^{1}$. Fix $\lambda>-\omega$. By assumption
(\ref{Axiom}), the fact that smooth compactly-supported tensor fields are
dense in $W_{2,0}^{1}$ implies that for all $u\in W_{2,0}^{1}$, one can
integrate by parts to get%
\[
\operatorname{Re}\left(  \beta_{\lambda}(u,u)\right)  \geq\left\Vert \nabla
u\right\Vert _{0}^{2}+(\lambda+\omega)\left\Vert u\right\Vert _{0}^{2}\geq
\min\{1,\lambda+\omega\}\left\Vert u\right\Vert _{1}^{2}.
\]
By the Lax--Milgram Theorem, for each $f\in L_{2,0}^{2}$, there exists a
unique $u\in W_{2,0}^{1}$ solving $\beta_{\lambda}(u,v)=\left(  f,\bar
{v}\right)  $ for all $v\in W_{2,0}^{1}$ and satisfying the estimate%
\begin{equation}
\left\Vert u\right\Vert _{1}\leq\max\{1,(\lambda+\omega)^{-1}\}\left\Vert
f\right\Vert _{0}. \label{LME}%
\end{equation}
It follows that $(\lambda+\omega)\left\Vert u\right\Vert _{0}^{2}%
\leq\operatorname{Re}\left(  \beta_{\lambda}(u,u)\right)  =\operatorname{Re}%
\left(  f,\bar{u}\right)  \leq\left\Vert f\right\Vert _{0}\left\Vert
u\right\Vert _{0}$, which implies that%
\[
\left\Vert u\right\Vert _{0}\leq\frac{1}{\lambda+\omega}\left\Vert
f\right\Vert _{0},
\]
and hence that%
\begin{equation}
\left\Vert \nabla u\right\Vert _{0}^{2}+(\lambda+\omega)\left\Vert
u\right\Vert _{0}^{2}\leq\operatorname{Re}(\left(  f,\bar{u}\right)
)\leq\frac{1}{\lambda+\omega}\left\Vert f\right\Vert _{0}^{2}.
\label{UpperBound}%
\end{equation}
Since the left-hand side of (\ref{UpperBound}) is an upper bound for
$\left\Vert u\right\Vert _{1}^{2}$ when $\lambda+\omega\geq1$, the conclusion
follows from (\ref{LME}) and (\ref{UpperBound}).
\end{proof}

Recall the linear operator $L$ defined by (\ref{L}). By Remark~\ref{BetterL},
$L$ can be written in the form%
\[
Lh=\Delta h-\nabla_{X}h+Zh,
\]
where $Z$ is an unbounded $(2,2)$-tensor field. For each soliton under
consideration, our calculations in Section \ref{Linear} show that $\delta X$
is constant. Moreover, by Propositions~\ref{Stable3}, \ref{StableSol}, and
\ref{Stable4}, there exists $\omega>0$ depending only on the soliton in
question such that
\begin{equation}
\left\langle \frac{1}{2}(\delta X)v-Zv,v\right\rangle \geq\omega\left\vert
v\right\vert ^{2} \label{GoodTerm}%
\end{equation}
holds pointwise for all real compactly supported $(2,0)$-tensor fields $v$.

Let $r=\sqrt{\sum_{\ell=1}^{n}x_{\ell}^{2}}$ denote the Euclidean distance
from the origin in $\mathbb{R}^{n}\approx\mathcal{M}^{n}$. By formulas
(\ref{X3}), (\ref{XSol}), and (\ref{X4}), there exists $C>0$ depending only on
the Lie group in question such that%
\begin{equation}
dr(X)\geq-Cr. \label{BadTerm}%
\end{equation}
For each $k=1,2,\ldots$, let $\gamma_{k}$ be a translated smooth bump function
such that $\gamma_{k}(s)=1$ for $s<k-1+\varepsilon$, $\gamma_{k}(s)=0$ for
$s>k-\varepsilon$, and $-2\leq\gamma^{\prime}(s)\leq0$ for all $s$. Define the
vector fields
\begin{equation}
\xi_{k}(x):=\gamma_{k}(r(x))X(x), \label{xi}%
\end{equation}
noting that $\xi_{k}$ is supported in the Euclidean ball of radius $k$. Define%
\begin{equation}
\zeta_{k}(x):=\left\{
\begin{array}
[c]{cl}%
Z(x) & \text{if }r(x)<k-1\\
\gamma_{k}(r(x))Z(x)-(Ck+\omega) & \text{if }r(x)\geq k-1
\end{array}
\right.  , \label{zeta}%
\end{equation}
where $C$ is given by (\ref{BadTerm}). Note that $\zeta_{k}$ is smooth except
on the hypersurface $r=k-1$.

\begin{lemma}
\label{Approximation}Assume (\ref{GoodTerm}) and (\ref{BadTerm}) hold. Fix any
$\lambda>-\omega$ and let $f\in L_{2,0}^{2}$ be given. Then for each integer
$k\geq1$, there exists a unique $u_{k}\in W_{2,0}^{1}$ solving%
\begin{equation}
\lambda\left(  u_{k},\bar{v}\right)  +\left(  \nabla u_{k},\nabla\bar
{v}\right)  +\left(  \nabla_{\xi_{k}}u_{k},\bar{v}\right)  -\left(  \zeta
_{k}u_{k},\bar{v}\right)  =\left(  f,\bar{v}\right)
\label{ApproximateSolution}%
\end{equation}
for all $v\in W_{2,0}^{1}$, where $\xi_{k}$ and $\zeta_{k}$ are defined in
(\ref{xi}) and (\ref{zeta}), respectively. Furthermore, all $u_{k}$ are
bounded uniformly by%
\[
\left\Vert u_{k}\right\Vert _{1}^{2}\leq\frac{1}{\lambda+\omega}\left\Vert
f\right\Vert _{0}^{2}.
\]

\end{lemma}

\begin{proof}
If $r<k-1$, then (\ref{GoodTerm}) gives the pointwise inequality
\[
\left\langle \frac{1}{2}(\delta\xi_{k})v-\zeta_{k}v,v\right\rangle
=\left\langle \frac{1}{2}(\delta X)v-Zv,v\right\rangle \geq\omega\left\vert
v\right\vert ^{2}%
\]
for all real $v$. If $r\geq k-1$, then we have%
\[
\left\langle \frac{1}{2}(\delta\xi_{k})v-\zeta_{k}v,v\right\rangle \geq
\gamma_{k}(r)\left\langle \frac{1}{2}(\delta X)v-Zv,v\right\rangle
+\omega\left\vert v\right\vert ^{2}\geq\omega\left\vert v\right\vert ^{2}%
\]
for all real $v$ by (\ref{GoodTerm}), (\ref{BadTerm}), the fact that $\xi_{k}$
is supported in $r<k$, and our choice of $\zeta_{k}$. Applying Lemma~\ref{LM}
with $\xi=\xi_{k}$ and $\zeta=\zeta_{k}$ thus yields the result.
\end{proof}

Now let $\mathbb{L}$ denote the complexification of the linear operator $L$
defined by (\ref{L}). $\mathbb{L}$ is defined for all $u=v+iw$ (with $v,w$
real) by%
\begin{equation}
\mathbb{L}u:=Lv+iLw.
\end{equation}
\emph{A priori, }$\mathbb{L}$ is only defined on $C_{0}^{\infty}(T_{2}%
^{0}\mathcal{M}^{n};\mathbb{C})$.

\begin{lemma}
\label{HY-input}Assume (\ref{GoodTerm}) and (\ref{BadTerm}) hold. Fix any
$\lambda>-\omega$. Then for each $f\in L_{2,0}^{2}$, there exists a weak
solution $u\in W_{2,0}^{1}$ of
\[
(\lambda\mathbb{I}-\mathbb{L})u=f
\]
such that%
\[
\left\Vert u\right\Vert _{1}^{2}\leq\frac{1}{\lambda+\omega}\left\Vert
f\right\Vert _{0}^{2}.
\]

\end{lemma}

\begin{proof}
Let $\{u_{k}\}_{k\geq1}$ be the sequence given by Lemma~\ref{Approximation}.
Since $u_{k}$ is bounded uniformly in $W_{2,0}^{1}$, a subsequence
$\{u_{k_{j}}\}_{j\in\mathbb{N}}$ converges weakly to some $u\in W_{2,0}^{1}$.
Thus the $(\lambda\mathbb{I}-\mathbb{L})u_{k_{j}}$ converge to $(\lambda
\mathbb{I}-\mathbb{L})u$ as distributions. This means that for any smooth
compactly-supported $(2,0)$-tensor field $v$, one has%
\[
\lim_{j\rightarrow\infty}\left(  (\lambda\mathbb{I}-\mathbb{L})u_{k_{j}}%
,\bar{v}\right)  =\left(  (\lambda\mathbb{I}-\mathbb{L})u,\bar{v}\right)  ,
\]
where we make sense of the Laplacian term by integration by parts.

Recall that $\xi_{k}=X$ and $\zeta_{k}=Z$ wherever $r<k-1$. Because the
support $\operatorname*{Sprt}(v)$ of $v$ is compact, there exists $J$
depending on $v$ such that
\begin{align*}
\left(  (\lambda\mathbb{I}-\mathbb{L})u_{k_{j}},\bar{v}\right)   &
=\lambda\left(  u_{k_{j}},\bar{v}\right)  +\left(  \nabla u_{k_{j}},\nabla
\bar{v}\right)  +\left(  \nabla_{X}u_{k_{j}},\bar{v}\right)  -\left(
Zu_{k_{j}},\bar{v}\right) \\
&  =\lambda\left(  u_{k_{j}},\bar{v}\right)  +\left(  \nabla u_{k_{j}}%
,\nabla\bar{v}\right)  +\left(  \nabla_{\xi_{k_{j}}}u_{k_{j}},\bar{v}\right)
-\left(  \zeta_{k_{j}}u_{k_{j}},\bar{v}\right) \\
&  =\left(  f,\bar{v}\right)
\end{align*}
for all $j\geq J$. The first equality above follows directly from the
definition of $\mathbb{L}$. The second equality holds because $\xi_{k_{j}%
}|_{\operatorname*{Sprt}(v)}=X|_{\operatorname*{Sprt}(v)}$ and $\zeta_{k_{j}%
}|_{\operatorname*{Sprt}(v)}=Z|_{\operatorname*{Sprt}(v)}$ for all $j\geq J$,
by our choice of $J$. The last equality is (\ref{ApproximateSolution}).
Combining these observations, we see that $\left(  (\lambda\mathbb{I}%
-\mathbb{L})u,\bar{v}\right)  =\left(  f,\bar{v}\right)  $ for all compactly
supported $v$. By density of smooth compactly-supported tensor fields, this is
possible only if $(\lambda\mathbb{I}-\mathbb{L})u=f$. Since weak convergence
$u_{k_{j}}\rightharpoonup u$ implies that
\[
\left\Vert u\right\Vert _{1}\leq\liminf_{j\rightarrow\infty}\left\Vert
u_{k_{j}}\right\Vert _{1},
\]
the result follows.
\end{proof}

Now for any $\lambda>-\omega$, we may define the domain of $\mathbb{L}$ by%
\begin{equation}
\mathrm{D}_{\mathbb{L}}=(\lambda\mathbb{I}-\mathbb{L})^{-1}(L_{2,0}^{2}).
\end{equation}
By the resolvent identity $(\mu\mathbb{I}-\mathbb{L})^{-1}=(\lambda
\mathbb{I}-\mathbb{L})^{-1}+(\lambda-\mu)(\lambda\mathbb{I}-\mathbb{L}%
)^{-1}(\mu\mathbb{I}-\mathbb{L})^{-1}$, $\mathrm{D}_{\mathbb{L}}$ is well
defined, i.e.~independent of $\lambda>-\omega$. Because it contains
$C_{0}^{\infty}(T_{2}^{0}\mathcal{M}^{n};\mathbb{C})$, the domain
$\mathrm{D}_{\mathbb{L}}$ is dense in each $W_{2,0}^{k}$. The following is
then an immediate consequence of Lemma~\ref{HY-input} and the Hille--Yosida
Theorem. (See \cite[Theorem~3.1]{Pazy83}.)

\begin{theorem}
\label{StronglyContinuous}$\mathbb{L}$ generates a unique $C_{0}$ semigroup
with domain $\mathrm{D}_{\mathbb{L}}$. In particular, for every $f\in
L_{2,0}^{2}$ and $\lambda\in\mathbb{C}$ with $\operatorname{Re}(\lambda
)>-\omega$, there exists $u\in W_{2,0}^{1}$ solving $u=(\lambda\mathbb{I}%
-\mathbb{L})^{-1}f$ and satisfying%
\[
\left\Vert u\right\Vert _{1}^{2}\leq\frac{1}{\operatorname{Re}(\lambda
)+\omega}\left\Vert f\right\Vert _{0}^{2}.
\]

\end{theorem}

\section{Directions for further research\label{Open}}

This note is a contribution to a challenging program to investigate analytic
aspects of large-time behavior of immortal (especially Type-III) Ricci flow
solutions, and in particular to investigate Conjectures~\ref{Old} and
\ref{New}. One aspect of this program may be divided into three main parts,
all of which are undergoing active development:

\begin{enumerate}
\item Catalog examples of expanding homogeneous Ricci solitons as models of
Type-III behavior.

\item Determine asymptotic stability of these homogeneous Ricci solitons.

\item Deduce analytic properties of their pseudosoliton quotients.
\end{enumerate}

Contributions to Part~1 have been made by Lauret \cite{Lauret01}, Lauret--Will
\cite{LW06}, Baird--Danielo \cite{BD06}, and Lott \cite{Lott06}, among others.
Results of Glickenstein \cite{Glickenstein03} and Lott \cite{Lott06}
contribute to Part~3, as do the theories of Riemannian groupoids and
megafolds. (See \cite{Haefliger01} and \cite{PT99}, for example.) Part~2 is
still largely open. We hope this note stimulates progress toward its resolution.


\begin{thebibliography}{99}                                                                                               %


\bibitem {AK1}\textbf{Angenent, Sigurd B.; Knopf, Dan.} An example of
neckpinching for Ricci flow on $S^{n+1}$. \emph{Math.~Res.~Lett. }\textbf{11}
(2004), no.~4, 493--518.

\bibitem {AK2}\textbf{Angenent, Sigurd B.; Knopf, Dan.} Precise asymptotics
for the Ricci flow neckpinch. \texttt{arXiv:math.DG/0511247}.

\bibitem {BD06}\textbf{Baird, Paul; Danielo, Laurent.} Three-dimensional Ricci
solitons which project to surfaces. \texttt{arXiv:\allowbreak math.DG/0510313}.

\bibitem {BTV95}\textbf{Berndt, J\"{u}rgen; Tricerri, Franco; Vanhecke,
Lieven.} \emph{Generalized Heisenberg groups and Damek-Ricci harmonic spaces.
}Lecture Notes in Mathematics, 1598. Springer-Verlag, Berlin, 1995.

\bibitem {Biquard00}\textbf{Biquard, Olivier.} M\'{e}triques d'Einstein
asymptotiquement sym\'{e}triques. \emph{Ast\'{e}risque } No.~265 (2000).

\bibitem {CFG92}\textbf{Cheeger, Jeff; Fukaya, Kenji; Gromov, Mikhael.}
Nilpotent structures and invariant metrics on collapsed manifolds.
\emph{J.~Amer.~Math.~Soc.~}\textbf{5} (1992), no.~2, 327--372.

\bibitem {CG86}\textbf{Cheeger, Jeff; Gromov, Mikhael.} Collapsing Riemannian
manifolds while keeping their curvature bounded.~I. \emph{J.~Differential
Geom.~}\textbf{23} (1986), no.~3, 309--346.

\bibitem {CG90}\textbf{Cheeger, Jeff; Gromov, Mikhael}. Collapsing Riemannian
manifolds while keeping their curvature bounded.~II. \emph{J.~Differential
Geom.~}\textbf{32} (1990), no.~1, 269--298.

\bibitem {CF04}\textbf{Cupini, Giovanni; Fornaro, Simona.} Maximal regularity
in $L^{p}(\mathbb{R}^{N})$ for a class of elliptic operators with unbounded
coefficients. \emph{Differential Integral Equations }\textbf{17} (2004),
no.~3-4, 259--296.

\bibitem {CV87}\textbf{Cannarsa, Piermarco; Vespri, Vincenzo.} Generation of
analytic semigroups by elliptic operators with unbounded coefficients.
\emph{SIAM J.~Math.~Anal.~}\textbf{18} (1987), no.~3, 857--872.

\bibitem {CZ06}\textbf{Cao, Huai-Dong; Zhu, Xi-Ping.} A Complete Proof of the
Poincar\'{e} and Geometrization Conjectures - application of the
Hamilton-Perelman theory of the Ricci flow. \emph{Asian J.~Math.~}\textbf{10}
(2006), no.~2, 165--492.

\bibitem {DPL88}\textbf{Da Prato, G.; Lunardi, A.} Stability, instability and
center manifold theorem for fully nonlinear autonomous parabolic equations in
Banach space. \emph{Arch.~Rational Mech.~Anal.~}\textbf{101} (1988), no.~2, 115--141.

\bibitem {DT83}\textbf{DeTurck, Dennis M.} Deforming metrics in the direction
of their Ricci tensors. \emph{J.\ Differential Geom.\ }\textbf{18} (1983),
no.\ 1, 157--162.

\bibitem {DT03}\textbf{DeTurck, Dennis M.} Deforming metrics in the direction
of their Ricci tensors, improved version. \emph{Collected Papers on Ricci
Flow. }Edited by H.-D.~Cao, B.~Chow, S.-C.~Chu, and S.-T.~Yau. Internat.
Press, Somerville, MA, 2003.

\bibitem {EF98}\textbf{Eichhorn, J\"{u}rgen; Fricke, Jan.} The module
structure theorem for Sobolev spaces on open manifolds.
Math.~Nachr.~\textbf{194} (1998), 35--47.

\bibitem {Fukaya87}\textbf{Fukaya, Kenji.} Collapsing Riemannian manifolds to
ones of lower dimensions. \emph{J.~Differential Geom.~}\textbf{25} (1987),
no.~1, 139--156.

\bibitem {Fukaya89}\textbf{Fukaya, Kenji.} Collapsing Riemannian manifolds to
ones with lower dimension.~II. \emph{J.~Math.~Soc.~Japan }\textbf{41} (1989),
no.~2, 333--356.

\bibitem {Glickenstein03}\textbf{Glickenstein, David.} Precompactness of
solutions to the Ricci flow in the absence of injectivity radius estimates.
\emph{Geom.~Topol.~}\textbf{7} (2003), 487--510.

\bibitem {GIK02}\textbf{Guenther, Christine; Isenberg, James; Knopf, Dan.}
Stability of the Ricci flow at Ricci-flat metrics. \emph{Comm.~Anal.~Geom.~}%
\textbf{10} (2002), no.~4, 741--777.

\bibitem {GMV02}\textbf{Gozzi, Fausto; Monte, Roberto; Vespri, Vincenzo.}
Generation of analytic semigroups and domain characterization for degenerate
elliptic operators with unbounded coefficients arising in financial
mathematics.~I. \emph{Differential Integral Equations }\textbf{15} (2002),
no.~9, 1085--1128.

\bibitem {Haefliger01}\textbf{Haefliger, Andr\'{e}.} Groupoids and foliations.
\emph{Groupoids in analysis, geometry, and physics (Boulder, CO, 1999),
}83--100, Contemp.~Math., \textbf{282}, Amer.~Math.~Soc., Providence, RI, 2001.

\bibitem {Hamilton95}\textbf{Hamilton, Richard S.} The formation of
singularities in the Ricci flow.\ \emph{Surveys in differential geometry,
Vol.~II }(Cambridge, MA, 1993), 7--136, Internat. Press, Cambridge, MA, 1995.

\bibitem {Hamilton95a}\textbf{Hamilton, Richard S.} A compactness property for
solutions of the Ricci flow. \emph{Amer.~J.~Math.~}\textbf{117} (1995), no.~3, 545--572.

\bibitem {Hamilton99}\textbf{Hamilton, Richard S.} Non-singular solutions of
the Ricci flow on three-manifolds. \emph{Comm.~Anal.~Geom.~}\textbf{7} (1999),
no.~4, 695--729.

\bibitem {HI93}\textbf{Hamilton, Richard; Isenberg, James.} Quasi-convergence
of Ricci flow for a class of metrics. \emph{Comm.~Anal.~Geom.~}\textbf{1}
(1993), no.~3-4, 543--559.

\bibitem {Heber98}\textbf{Heber, Jens.} Noncompact homogeneous Einstein
spaces. \emph{Invent.~Math.~}\textbf{133} (1998), no.~2, 279--352.

\bibitem {IJ92}\textbf{Isenberg, James; Jackson, Martin.} Ricci flow of
locally homogeneous geometries on closed manifolds. \emph{J.~Differential
Geom.~}\textbf{35} (1992), no.~3, 723--741.

\bibitem {IJL05}\textbf{Isenberg, James; Jackson, Martin; Lu, Peng.} Ricci
flow on locally homogeneous closed 4-manifolds. \texttt{arXiv:math.DG/0502170}.

\bibitem {Ivey93}\textbf{Ivey, Thomas.} Ricci solitons on compact
three-manifolds. \emph{Differential Geom.~Appl.~}\textbf{3} (1993), no.~4, 301--307.

\bibitem {KL}\textbf{Kleiner, Bruce; Lott, John.} Notes on Perelman's papers.
(25/05/2006 version) \texttt{arXiv:math.DG/0605667}.

\bibitem {Knopf00}\textbf{Knopf, Dan.} Quasi-convergence of the Ricci flow.
\emph{Comm.\ Anal.\ Geom.~}\textbf{8} (2000), no.~2, 375--391.

\bibitem {KM01}\textbf{Knopf, Dan; McLeod, Kevin}. Quasi-convergence of model
geometries under the Ricci flow. \emph{Comm.\ Anal.\ Geom.}~\textbf{9} (2001),
no.~4, 879--919.

\bibitem {Koiso79}\textbf{Koiso, Norihito.} On the second derivative of the
total scalar curvature. \emph{Osaka J.~Math. }\textbf{16} (1979), no.~2, 413--421.

\bibitem {Lauret01}\textbf{Lauret, Jorge.} Ricci soliton homogeneous
nilmanifolds. \emph{Math.~Ann.~}\textbf{319} (2001), no.~4, 715--733.

\bibitem {LW06}\textbf{Lauret, Jorge.; Will, Cynthia. }Einstein solvmanifolds:
existence and non-existence questions. \texttt{arXiv:math.DG/0602502}.

\bibitem {LMP05}\textbf{Lunardi, A.; Metafune, G.; Pallara, D.} Dirichlet
boundary conditions for elliptic operators with unbounded drift.
\emph{Proc.~Amer.~Math.~Soc.~}\textbf{133} (2005), no.~9, 2625--2635.

\bibitem {Lott06}\textbf{Lott, John.} On the long-time behavior of type-III
Ricci flow solutions. \texttt{arXiv:\allowbreak math.DG/\allowbreak0509639}.

\bibitem {Milnor76}\textbf{Milnor, John.} Curvatures of left invariant metrics
on Lie groups. \emph{Advances in Math.~}\textbf{21} (1976), no.~3, 293--329.

\bibitem {MPSR05}\textbf{Metafune, G.; Pr\"{u}ss, J.; Schnaubelt, R.; Rhandi.
A.} $L^{p}$-regularity for elliptic operators with unbounded coefficients.
\emph{Adv.~Differential Equations }\textbf{10} (2005), no.~10, 1131--1164.

\bibitem {Pazy83}\textbf{Pazy, Ammon.} \emph{Semigroups of linear operators
and applications to partial differential equations. }Applied Mathematical
Sciences, \textbf{44}. Springer-Verlag, New York, 1983.

\bibitem {Perelman1}\textbf{Perelman, Grisha.} The entropy formula for the
Ricci flow and its geometric applications. \texttt{arXiv:math.DG/0211159}.

\bibitem {Perelman2}\textbf{Perelman, Grisha.} Ricci flow with surgery on
three-manifolds. \texttt{arXiv:math.DG/0303109}.

\bibitem {PT99}\textbf{Petrunin, A.; Tuschmann, W.} Diffeomorphism finiteness,
positive pinching, and second homotopy. \emph{Geom.~Funct.~Anal.~}\textbf{9}
(1999), no.~4, 736--774.

\bibitem {Rabier05}\textbf{Rabier, Patrick J.} Elliptic problems on
$\mathbb{R}^{N}$ with unbounded coefficients in classical Sobolev spaces.
\emph{Math.~Z.~}\textbf{249} (2005), no.~1, 1--30.

\bibitem {Sesum04}\textbf{\v{S}e\v{s}um, Nata\v{s}a}. Linear and dynamical
stability of Ricci flat metrics. \texttt{arXiv:math.DG\allowbreak/0410062}.

\bibitem {SY00}\textbf{Shioya, Takashi; Yamaguchi, Takao.} Collapsing
three-manifolds under a lower curvature bound. \emph{J.~Differential
Geom.~}\textbf{56} (2000), no.~1, 1--66.

\bibitem {SY05}\textbf{Shioya, Takashi; Yamaguchi, Takao.} Volume collapsed
three-manifolds with a lower curvature bound. \emph{Math.~Ann.~}\textbf{333}
(2005), no.~1, 131--155.

\bibitem {Schueth2004}\textbf{Schueth, Dorothee.} On the `standard' condition
for noncompact homogeneous Einstein spaces. \emph{Geom.~Dedicata }\textbf{105}
(2004), 77--83.
\end{thebibliography}
\end{document}